\documentclass{article}

\usepackage[english]{babel}

\usepackage[letterpaper,top=2cm,bottom=2cm,left=3cm,right=3cm,marginparwidth=1.75cm]{geometry}

\usepackage{amsmath,amsthm,amssymb,mathrsfs,bm}
\usepackage{graphicx}
\usepackage[colorlinks=true, allcolors=blue]{hyperref}
\usepackage{doi}
\numberwithin{equation}{section}

\newtheorem{theorem}{\bf{Theorem}}[section]
\newtheorem{lemma}{\bf {Lemma}}[section]
\newtheorem{define}{\bf{Definition}}[section]
\newtheorem{coro}{\bf{Corollary}}[section]

\newtheorem{remark}{\bf{Remark}}[section]

\newcommand{\be}{\begin{equation}}
\newcommand{\ee}{\end{equation}}
\newcommand\bes{\begin{eqnarray}}
\newcommand\ees{\end{eqnarray}}
\newcommand{\bess}{\begin{eqnarray*}}
\newcommand{\eess}{\end{eqnarray*}}
\newcommand{\mbE}{\hat{\mathbb{E}}}
\newcommand{\mbe}{\hat{\mathcal{E}}}
\newcommand{\tE}{\tilde{\mathbb{E}}}
\newcommand{\V}{\mathbb{V}}
\newcommand{\mv}{\mathcal{V}}
\newcommand{\sles}{(\Omega,\mathcal{H},\mbE)}
\newcommand{\usigma}{\overline{\sigma}}
\newcommand{\lsigma}{\underline{\sigma}}

\newcommand{\Xnk}{X_{n,k}}
\newcommand{\ld}{\underline{\delta}}
\newcommand{\ud}{\overline{\delta}}

\title{Central Limit Theorem for $m$-dependent random variables under sub-linear expectations}
\author{Wang-Yun Gu,\enspace Li-Xin Zhang}
\date{}

\begin{document}
\maketitle

\begin{abstract}
\par $M$-dependence is a commonly used assumption in the study of dependent sequences. In this paper, central limit theorems for $m$-dependent random variables under the sub-linear expectations are established based mainly on the conditions of Zhang\cite{zhang23}. They can be regarded as the extension of independent Lindeberg central limit theorem and for proving this, Rosenthal's inequality for $m$-dependent random variables is obtained. In particular, we extend the results in Li\cite{li15} and establish the central limit theorem for $m$-dependent stationary sequence.\\
\par \bf{Keywords:}\rm\quad central limit theorem, $m$-dependence, sub-linear expectation, capacity.
\end{abstract}

\section{Introduction and notations}
In the framework of sub-linear expectations introduced by Peng\cite{Peng19}, the expectation function is a sub-linear function instead of linear function, thus both the expectation and the related probability, called capacity, are non-additive, which are useful for studying models with uncertainty. Central limit theorem is of great significance in classical probability theory and is widely used in statistics, finance and many other fields. For the central limit theorems under the sub-linear expectations, Peng\cite{Peng08,Peng19b} and Krylov\cite{Krylov20} obtained the central limit theorem for independent and identically distributed(i.i.d.) random variables. Zhang\cite{zhang20} established the Lindeberg central limit theorem and functional central limit theorem for independent but not necessarily identically distributed one-dimensional random variables as well as martingale-like sequence. As a promotion, Zhang\cite{zhang23} proved that the Lindeberg central limit theorems still hold for multi-dimensional random vectors. 
\par However, independence of random variables is not usually satisfied in the real world, m-dependence is a weak and reasonable condition. Though central limit theorem for martingale-like sequence can deal with dependence case, we still hope to establish a specific theorem for $m$-dependent random variables with conditions without conditional expecatations. In the classical probability theory, Hoeffding and Robbins\cite{HR48}, Diananda\cite{Dia55} and Orey\cite{Orey58} proved the central limit theorem for $m$-dependent random variables. Furthermore, Berk\cite{Berk73} and Romano and Wolf\cite{RW00} established the theorems for the case $m$ increasing to infinity with $n$ under some additional conditions. Li\cite{li15} obtained the central limit theorem for $m$-dependent random variables under the sub-linear expectations based on the theorem of Peng\cite{Peng08}. In this paper, we derive the central limit theorem for one-dimensional $m$-dependent random variables under the Lindeberg condition, which improves the results in Li\cite{li15}. As a corollary, a special case for a stationary sequence of $m$-dependent random variables without mean uncertianty is considered and only simple conditions are needed for this case. We also establish the $m$-dependent central limit theorem with conditions for truncated random variables and a condition of capacity, which is weaker than Lindeberg condition. To deal with the terms of product of independent random variables with mean uncertainty and verify the convergence for functions with polynomial growth, we also prove the Rosenthal's inequality for $m$-dependent random variables under the sub-linear expectations.
\par In the rest of this section, we state some notations about sub-linear expectation. In Section 2, we cite the one-dimensional central limit theorem of Zhang\cite{zhang23} and prove Rosenthal's inequality and central limit theorems for both array and sequence of one-dimensional $m$-dependent random variables. We also derive the central limit theorem for $m$-dependent stationary sequence in this section. We give central limit theorems for truncated independent and $m$-dependent random variables in Section 3.
\par We use the framework and notations of Peng\cite{Peng08,Peng09,Peng19}.If one is familiar with these notations, he or she can skip this section. Let $(\Omega,\mathcal{F})$ be a given measurable space and let $\mathcal{H}$ be a linear space of real functions defined on  $(\Omega,\mathcal{F})$ such that if $X_1,\cdots, X_n\in\mathcal{H}$, then $\varphi(X_1,\cdots,X_n)\in\mathcal{H}$ for each $\varphi\in C_{l,Lip}(\mathbb{R}^n)$, where $C_{l,Lip}(\mathbb{R}^n)$ denotes the linear space of local Lipschitz functions $\varphi$ satisfying
\bess
\vert \varphi(\bm{x})-\varphi(\bm{y})\vert\leq C(1+\vert \bm{x}\vert^m+\vert\bm{y}\vert^m)\vert\bm{x}-\bm{y}\vert,\quad\forall\bm{x},\bm{y}\in\mathbb{R}^n,
\eess
for some $C>0, m\in\mathbb{N}$ depending on $\varphi$.\\
$\mathcal{H}$ is considered as a space of "random variables". We also denote $C_{b,Lip}(\mathbb{R}^n)$ the space of bounded Lipschitz functions. In this case, we denote $X\in\mathcal{H}$. 
\begin{define}
A sub-linear expectation $\mbE$ on $\mathcal{H}$ is a function $\mbE:\mathcal{H}\rightarrow\bar{\mathbb{R}}$ satisfying the following properties: for all $X,Y\in\mathcal{H}$, we have
\begin{itemize}
\item[(a)] Monotonicity: If $X\geq Y$, then $\mbE[X]\geq\mbE[Y]$;
\item[(b)] Constant preserving: $\mbE[c]=c$;
\item[(c)] Sub-additivity: $\mbE[X+Y]\leq\mbe[X]+\mbE[Y]$ whenever $\mbE[X]+\mbE[Y]$ is not of the form $+\infty-\infty$ or $-\infty+\infty$;
\item[(d)] Positive homogeneity: $\mbE[\lambda X]=\lambda\mbE[X]$ for $\lambda\geq0$.
\end{itemize}
Here, $\bar{\mathbb{R}}=[-\infty,\infty]$, $0\cdot\infty$ is defined to be 0. The triple $\sles$ is called a sub-linear expectation space. Give a sub-linear expectation $\mbE$, let us denote the conjugate expectation $\mbe$ of $\mbE$ by
\bess
\mbe[X]:=-\mbE[-X],\quad\forall X\in\mathcal{H}.
\eess
\end{define}
\par From the definition, it is easily shown that $\mbe[X]\leq\mbE[X],\enspace \mbE[X+c]=\mbE[X]+c$, and $\mbE[X-Y]\geq\mbE[X]-\mbE[Y]$ for all $X,Y\in\mathcal{H}$ with $\mbE[Y]$ being finite. We also call $\mbE[X]$ and $\mbe[X]$ the upper-expecation and lower-expectation of $X$, respectively.

\begin{define}
\begin{itemize}
\item[(i)] (Identical distribution) Let $\bm{X}_1$ and $\bm{X}_2$ be two n-dimensional random vectors, respectively, defined in sub-linear expectation spaces $(\Omega_1,\mathcal{H}_1,\mbE_1)$ and $(\Omega_2,\mathcal{H}_2,\mbE_2)$. They are called identically distributed, denoted by $\bm{X}_1\overset{d}{=}\bm{X}_2$, if
\bess
	\mbE_1[\varphi(\bm{X}_1)]=\mbE_2[\varphi(\bm{X}_2)],\quad\forall\varphi\in C_{l,Lip}(\mathbb{R}^n).
\eess
\par A sequence $\{X_n;n\geq 1\}$ of random variables is said to be identically distributed if $X_i\overset{d}{=}X_1$ for each $i\geq 1$.
\item[(ii)] (Independence) In a sub-linear expectation space $\sles$, a random vector $\bm{Y}=(Y_1,\cdots,Y_n),Y_i\in\mathcal{H}$ is said to be independent of another random vector $\bm{X}=(X_1,\cdots,X_m),X_i\in\mathcal{H}$ under $\mbE$ if for each test function $\varphi\in C_{l,Lip}(\mathbb{R}^m\times\mathbb{R}^n)$ we have $\mbE[\varphi(\bm{X},\bm{Y})]=\mbE[\mbE[\varphi(\bm{x},\bm{Y})]\vert_{\bm{x}=\bm{X}}]$, whenever $\bar{\varphi}(\bm{x}):=\mbE[\vert\varphi(\bm{x},\bm{Y})\vert]<\infty$ for all $\bm{x}$ and $\mbE[\vert\bar{\varphi}(\bm{X})\vert]<\infty.$
\item[(iii)] (Independent random variables) A sequence of random variables(or random vectors) $\{X_n;n\geq 1\}$ is said to be independent if $X_{i+1}$ is independent of $(X_1,\cdots,X_i)$ for each $i\geq 1$.
\item[(iv)](m-dependence) A sequence of random variables(or random vectors) $\{X_n;n\geq 1\}$ is said to be $m$-dependent if there exists an integer $m$ such that for every $n$ and every $j\geq m+1$, $(X_{n+m+1},\cdots,X_{n+j})$ is independent of $(X_1,\cdots,X_n)$. In particular, if $m=0$, $\{X_n;n\geq 1\}$ is an independent sequence.
\item[(v)](stationary) A sequence of random variables(or random vectors) $\{X_n;n\geq 1\}$ is said to be stationary if for every positive integers $n$ and $p$, $(X_1,\cdots,X_n)\overset{d}{=}(X_{1+p},\cdots,X_{n+p})$.
\end{itemize}
\end{define}
It is easily seen that if $\{X_1,\cdots,X_n\}$ are independent, then $\mbE[\sum_{i=1}^nX_i]=\sum_{i=1}^n\mbE[X_i]$.
\par Next, we consider the capacities corresponding to the sub-linear expectations. Let $\mathcal{G}\subset\mathcal{F}$. A function $V:\mathcal{G}\rightarrow[0,1]$ is called a capacity if
$$
V(\emptyset)=0,\enspace V(\Omega)=1 \enspace and\enspace  V(A)\leq V(B)\enspace \forall A\subset B, A,B\in\mathcal{G}.
$$
It is called sub-additive if $V(A\cup B)\leq V(A)+V(B)$ for all $A,B\in\mathcal{G}$ with $A\cup B\in\mathcal{G}$.
\par Let $\sles$ be a sub-linear expectation space. We define $(\V,\mv)$ as a pair of capacities with the properties that
\bes
	\mbE[f]\leq\V(A)\leq\mbE[g]\quad if\enspace f\leq I_A\leq g,f,g\in\mathcal{H}\enspace and\enspace A\in\mathcal{F},\label{eq1.1}
\ees
$\V$ is sub-additive and $\mv(A):=1-\V(A^c),A\in\mathcal{F}$. It is obvious that
\bes
\mv(A\cup B)\leq\mv(A)+\V(B).
\ees
We call $\V$ and $\mv$ the upper and lower capacity, respectively. In general, we choose $(\V,\mv)$ as
\bes
\hat{\V}(A):=\inf\{\mbE[\xi]:I_A\leq\xi,\xi\in\mathcal{H}\},\hat{\mv}(A)=1-\hat{\V}(A^c), \enspace\forall A\in\mathcal{F}.
\ees If $\V$ on the sub-linear expectation space $\sles$ and $\tilde{\mathbb{V}}$ on the sub-linear expectation space $(\tilde{\Omega},\tilde{\mathcal{H}},\tilde{\mathbb{E}})$ are two capacities have the property (\ref{eq1.1}), then for any random variables $X\in\mathcal{H}$ and $\tilde{X}\in\tilde{\mathcal{H}}$ with $X\overset{d}{=}\tilde{X}$, we have
\bes
	\V(X\geq x+\epsilon)\leq\tilde{V}(\tilde{X}\geq x)\leq\V(X\geq x-\epsilon)\quad for\enspace all\enspace\epsilon>0\enspace and \enspace x.\label{eq1.2}
\ees
In fact, let $f\in C_{b,Lip}(\mathbb{R})$ such that $I\{y\geq x+\epsilon\}\leq f(y)\leq I\{y\geq x\}$. Then
\bess
	\V(X\geq x+\epsilon)\leq\mbE[f(X)]=\tilde{\mathbb{E}}[f(\tilde{X})]\leq\tilde{\V}(X\geq x),
\eess
and similar $\tilde{V}(\tilde{X}\geq x)\leq\V(X\geq x-\epsilon)$. It follows from (\ref{eq1.2}) that
\bess
	\V(X\geq x)=\tilde{\V}(X\geq x),\enspace \V(X> x)=\tilde{\V}(X> x)
\eess
for all but except countable many $x$. In this paper the events that we considered are almost of the type $\{X\geq x\}$ of $\{X>x\}$, so the choice of capacity will not influence our results.
\par Moreover, we recall the definitions of types of convergence.
\begin{define}
\begin{itemize}
\item[(i)] A sequence of $d$-dimensional random vectors $\{\bm{X}_n;n\geq 1\}$ defined on a sub-linear expectation space $\sles$ is said to converge in distribution(or converge in law) under $\mbE$ if for each $\varphi\in C_{b,Lip}(\mathbb{R}^n)$, the sequence $\{\mbE[\varphi(\bm{X_n})];n\geq 1\}$ converges.
\item[(ii)]  A sequence of $d$-dimensional random vectors $\{\bm{X}_n;n\geq 1\}$ defined on a sub-linear expectation space $\sles$ is said to converge in $\V$ if there exists a $\mathcal{F}$-measurable random vector $\bm{X}$ such that
\bess
	\V(|\bm{X}_n-\bm{X}|>\epsilon)\rightarrow 0 \enspace \forall\epsilon>0.
\eess
\end{itemize}
\end{define}

\par Finally, we state the notations of G-normal distribution. Let $\mathbb{S}(d)$ be the collection of all $d\times d$ symmetric matrices. A function $G:\mathbb{S}(d)\rightarrow\mathbb{R}$ is called a sub-linear function monotonic in $A\in\mathbb{S}(d)$ if for each $A,\bar{A}\in\mathbb{S}(d)$,
\begin{align*}
		\begin{split}
 			\left \{
 			\begin{array}{ll}
 				&G(A+\bar{A})\leq G(A)+G(\bar{A}),\\
    			&G(\lambda A)=\lambda G(A),\enspace\forall\lambda>0,\\
				&G(A)\geq G(\bar{A}),\enspace if \enspace A\geq\bar{A}.
			\end{array}
		\right.
		\end{split}
\end{align*}
Here $A\geq\bar{A}$ means that $A-\bar{A}$ is semi-positive definitive. $G$ is continuous if $|G(A)-G(\bar{A})|\rightarrow0$ when $\Vert A-\bar{A}\Vert_\infty\rightarrow0$, where $\Vert A-\bar{A}\Vert_\infty=\max_{i,j}|a_{ij}-\bar{a}_{ij}|$ for $A=(a_{ij})_{i,j=1}^d$ and $\bar{A}=(\bar{a}_{ij})_{i,j=1}^d$.
\begin{define}
(G-normal random variable) Let $G:\mathbb{S}(d)\rightarrow\mathbb{R}$ be a continuous sub-linear function monotonic in $A\in\mathbb{S}(d)$. A $d$-dimensional random vector $\bm{\xi}=(\xi_1,\cdots,\xi_d)$ in a sub-linear expectation space $(\tilde{\Omega},\tilde{\mathcal{H}},\tE)$ is called a G-normal distributed random variable(written as $\xi\sim N(0,G)$ under $\tE$), if for any $\varphi\in C_{l,Lip}(\mathbb{R}^d)$, the function $u(t,\bm{x})=\tE[\varphi(\bm{x}+\sqrt{t}\bm{\xi})]$($\bm{x}\in\mathbb{R}^d,t\geq0$) is the unique viscosity solution of the following heat equation:
\bess
	\partial_tu-\frac12G(D^2u)=0,\enspace u(0,\bm{x})=\varphi(\bm{x}),
\eess
where $Du=(\partial_{x_i}u,i=1,\cdots,d)$ and $D^2u=D(Du)=(\partial_{x_i,x_j}u)_{i,j=1}^d$.
\end{define}
That $\bm{\xi}$ is a G-normal distributed random vector is equivalent to that, if $\bar{\bm{\xi}}$ is an independent copy of $\bm{\xi}$,then
\bess
a\bm{\xi}+b\bar{\bm{\xi}}\overset{d}{=}\sqrt{a^2+b^2}\bm{\xi},\enspace\forall a,b\geq0,
\eess 
and $G(A)=\tE[\langle A\bm{\xi},\bm{\xi}\rangle]$(c.f. Definition 2.2.4 and Corollary 2.2.13 of Peng\cite{Peng19}), where $\langle\bm{x},\bm{y}\rangle$ is the scalar product of $\bm{x},\bm{y}$. When $d=1$, $G$ can be written as $G(\alpha)=\alpha^+\usigma^2-\alpha^-\lsigma^2$ and we write $\xi\sim N(0,[\lsigma^2,\usigma^2])$ if $\xi$ a G-normal distributed random variable.
\par Through this paper, for real numbers $x$ and $y$, we denote $x\vee y=\max\{x,y\},x\wedge y=\min\{x,y\},x^+=x\vee 0$ and $x^-=x\wedge 0$. For a random variable $X$, because $XI\{\vert X\vert\leq c\}$ may not be in $\mathcal{H}$, we will truncate it in the form $(-c)\vee X\wedge c$ denoted by $X^{(c)}$.
In the sequel, the constants $C,C_p$ and $C_{m,p}$ can represent different values from line to line.

\section{Central Limit Theorem for $m$-dependent random variables}
In this section, we prove the central limit theorem for arrays of $m$-dependent random variables. First, we introduce the Corollary 3.3 of Zhang\cite{zhang23} and prove the Rosenthal-type inequality for $m$-dependent random variables based on Theorem 2.1 of Zhang\cite{zhang16}.
\begin{lemma}
	Let $\{X_{n,k};k=1,\cdots,k_n\}$ be an array of independent random variables, $n=1,2,\cdots$. Denote $\usigma_{n,k}^2=\mbE[X_{n,k}^2],\lsigma_{n,k}^2=\mbe[X_{n,k}^2]$ and $B_n^2=\sum_{k=1}^{k_n}\usigma_{n,k}^2$. Suppose that the Lindeberg condition is satisfied:
\bes
	\frac{1}{B_n^2}\sum_{k=1}^{k_n}\mbE\left[(X_{n,k}^2-\epsilon B_n^2)^+\right]\rightarrow 0 \enspace\forall \epsilon>0,\label{Lin}
\ees
and further, there is a constant $r\in[0,1]$ such that 
\bes
\frac{\sum_{k=1}^m\lsigma_{n,k}^2}{\sum_{k=1}^m\usigma_{n,k}^2}\rightarrow r,\label{var}
\ees
as long as $k_n>m\rightarrow\infty$ with $\liminf\frac{\sum_{k=1}^m\usigma_{n,k}^2}{B_n^2}>0$, and
\bes
	\frac{\sum_{k=1}^{k_n}\{|\mbE[X_{n,k}]|+|\mbe[X_{n,k}]|\}}{B_n}\rightarrow 0.\label{mean}
\ees
Then for any continuous function $\varphi$ with $|\varphi(x)|\leq Cx^2$,
\bes
\lim_{n\rightarrow\infty}\mbE\left[\varphi\left(\frac{\sum_{k=1}^{k_n}X_{n,k}}{B_n}\right)\right]=\tE[\varphi(\xi)],
\ees
where $\xi\sim N(0,[r,1])$ under $\tE$.\label{lemma1d}
\end{lemma}
\begin{lemma}
	Suppose $\{X_k;k=1,\cdots,n\}$ is a sequence of $m$-dependent random variables in the sub-linear space $\sles$ . Denote $S_k=\sum_{i=1}^kX_i$, then for $p\geq2$,
\bes
	\mbE\left[\max_{k\leq n}|S_k|^p\right]\leq C_{m,p}\left\{\sum_{k=1}^n\mbE[|X_k|^p]+\left(\sum_{k=1}^n\mbE[|X_k|^2]\right)^{p/2}+\left(\sum_{k=1}^n\left[|\mbE[X_k]|+|\mbe[X_k]|\right]\right)^p\right\}.\label{R1}
\ees
Here $C_{m,p}$ is a positive constant depending on $m$ and $p$.
\label{lemma4.1}
\end{lemma}
\proof Denote $I_{j,k}=\{i\in\mathbb{N}:i\leq k,\enspace i\enspace mod\enspace (m+1)\equiv j\},j=0,\cdots,m$ and $S_{j,k}=\sum_{i\in I_{j,k}}X_{i}$. Note that
\bess
	\left|S_k\right|^p=\left|\sum_{j=0}^m S_{j,k}\right|^p\leq (m+1)^p\max_{0\leq j\leq m}|S_{j,k}|^p\leq (m+1)^p\sum_{j=0}^m|S_{j,k}|^p,
\eess
 we have
\bess
	\mbE\left[\max_{k\leq n}|S_k|^p\right]\leq C_{m,p}\mbE\left[\max_{k\leq n}\left(\sum_{j=1}^m|S_{j,k}|^p\right)\right]\leq C_{m,p}\sum_{j=0}^m\mbE\left[\max_{k\leq n}|S_{j,k}|^p\right].
\eess
Since $S_{j,k}$s are partial sums of  independent random variables and all the terms on the right side of Rosenthal's inequality for independent random variables are positive, we obtain (\ref{R1}) easily.\qedsymbol\\

With Rosenthal's inequality for $m$-dependent random variables and the construction of division by Orey\cite{Orey58}, we have the central limit theorem for $m$-dependent random variables.
\begin{theorem}
	Let $\{X_{n,k};k=1,\cdots,k_n\}$ be an array of  $m$-dependent random variables in the sub-linear expectation space $\sles$ with $\mbE[X_{n,k}^2]<\infty$. Denote $B_n^2=\mbE\left[\left(\sum_{k=1}^{k_n}X_{n,k}\right)^2\right]$. Assume that 
\bes
	&&\frac{1}{B_n^2}\sum_{k=1}^{k_n}\mbE[(X_{n,k}^2-\epsilon B_n^2)^+]\rightarrow0\enspace\forall\epsilon>0,\label{eq3.26}\\
	&&\frac{1}{B_n}\sum_{k=1}^{k_n}\left\{|\mbE[X_{n,k}]|+|\mbe[X_{n,k}]|\right\}\rightarrow0,\label{eq3.27}\\
	&&\frac{1}{B_n^2}\sum_{k=1}^{k_n}\mbE[X_{n,k}^2]=O(1),\label{eq3.28}
\ees
and there exist a constant $r\in[0,1]$ such that 
\bes
	\frac{\mbe\left[\left(\sum_{k=1}^{M}X_{n,k}\right)^2\right]}{\mbE\left[\left(\sum_{k=1}^{M}X_{n,k}\right)^2\right]}\rightarrow r,\label{eq3.29}
\ees
whenever $k_n>M\rightarrow\infty$ with $\liminf\frac{\mbE\left[\left(\sum_{k=1}^{M}X_{n,k}\right)^2\right]}{B_n^2}>0$.
Then for any continuous function $\varphi$ with $|\varphi(x)|\leq Cx^2$, we have
\bes
	\lim_{n\rightarrow\infty}\mbE\left[\varphi\left(\frac{\sum_{k=1}^{k_n}X_{n,k}}{B_n}\right)\right]=\tE[\varphi(\xi)],\label{eq3.30}
\ees
where $\xi\sim N(0,[r,1])$. Further, when $p>2$, (\ref{eq3.30}) holds for any continuous function $\varphi$ with $|\varphi(x)|\leq C|x|^p$ if (\ref{eq3.26}) is replaced by the condition that 
\bes
	\frac{1}{B_n^p}\sum_{k=1}^{k_n}\mbE[|X_{n,k}|^p]\rightarrow0.\label{eqpp}
\ees
\label{th3.3}
\end{theorem}
\proof We first suppose that $m=1$. It follows from (\ref{eq3.26}) that there exists a sequence of positive numbers $p_n'\uparrow\infty$ such that
\bess
	\frac{p_n'^2}{B_n^2}\sum_{k=1}^{k_n}\mbE[(\Xnk^2-\frac{B_n^2}{p_n'^2})^+]\rightarrow 0, 
\eess
let $p_n$ be the largest even integer satisfying $p_n\uparrow\infty$ and $p_n\leq \sqrt{p_n'}$, it is obvious that
\bess
	\frac{p_n^4}{B_n^2}\sum_{k=1}^{k_n}\mbE[(\Xnk^2-\frac{B_n^2}{p_n^4})^+]\rightarrow 0.
\eess
When $n$ is large enough, it follow that $\frac{\epsilon^2}{p_n^2}\geq\frac{1}{p_n^4}$, thus for any $\epsilon>0$,
\bess
	\frac{p_n^2}{B_n^2}\sum_{k=1}^{k_n}\mbE[(\Xnk^2-\frac{\epsilon^2B_n^2}{p_n^2})^+]\leq \epsilon^2\frac{p_n^4}{B_n^2}\sum_{k=1}^{k_n}\mbE[(\Xnk^2-\frac{B_n^2}{p_n^4})^+]\rightarrow 0.
\eess
In the sequel, we set $\Xnk=0$ if $k=0$ or $k>k_n$. For $k=1,\cdots,k_n$, denote 
\bess
\beta_{n,k}&=&\{\mbE[X_{n,k-1}^2]+\mbE[\Xnk^2]+\mbE[X_{n,k+1}^2]\}/B_n^2,\\
\ld_{n,k}&=&\{\mbe[\Xnk^2]+2\mbe[X_{n,k}X_{n,k-1}]+2\mbe[X_{n,k}X_{n,k+1}]\}/B_n^2,\\
\ud_{n,k}&=&\{\mbE[\Xnk^2]+2\mbE[X_{n,k}X_{n,k-1}]+2\mbE[X_{n,k}X_{n,k+1}]\}/B_n^2.
\eess
Set $P_0=\{0\},g(0)=0$ and we define $P_{i+1},g(i+1)$ recursively until $g(i)+p_n>k_n$ by
\bess
	P_{i+1}=\{k\in\mathbb{N}:g(i)+p_n/2<k\leq g(i)+p_n\},\enspace g(i+1)=\min_{j\in P_{i+1}}\{j:\beta_{n,j}=\min_{l\in P_{i+1}}\beta_{n,l}\}.
\eess
Suppose $P_i,i=1,g(i),i=1,\cdots,h_n-1$ are defined, then we have
\bess
	g(h-1)+p_n>k_n,\enspace g(i-1)+p_n/2<g(i)\leq g(i-1)+p_n,\enspace i=1,\cdots,h-1.
\eess
From the construction of $P$ and $g$, there exist positive integers $r_{i,s}$, $i=1,\cdots,h-1$, $s=1,\cdots,p_n/2$ such that
\bess
	(1)\enspace r_{i,s}\in P_i,\quad(2)\enspace r_{i,s}\neq r_{i,t}\enspace s\neq t,\quad (3)\enspace \sum_{i=1}^{h-1}\beta_{n,g(i)}\leq\sum_{i=1}^{h-1}\beta_{n,r_{i,s}},\forall s.
\eess
So
\bess
\frac{p_n}{2}\sum_{i=1}^{h-1}\beta_{n,g(i)}\leq\sum_{i=1}^{h-1}\sum_{s=1}^{p_n/2}\beta_{n,r_{i,s}}\leq\sum_{k=1}^{k_n}\beta_{n,k}\leq\frac{3}{B_n^2}\sum_{k=1}^{k_n}\mbE[X_{n,k}^2]=O(1).
\eess
It follows that
\bes
	\sum_{i=1}^{h-1}\beta_{n,g(i)}\rightarrow0,\label{eq3.37}
\ees
and
\bes
	\left|\sum_{i=1}^{h-1}\ld_{n,g(i)}\right|&\leq&\frac{1}{B_n^2}\sum_{i=1}^{h-1}\{\mbe[X_{n,g(i)}^2]+2|\mbe[X_{n,g(i)}X_{n,g(i)-1}]|+2|\mbe[X_{n,g(i)}X_{n,g(i)+1}]|\}\nonumber\\
&\leq&\frac{1}{B_n^2}\sum_{i=1}^{h-1}\{\mbE[X_{n,g(i)}^2]+2\mbE[|X_{n,g(i)}X_{n,g(i)-1}|]+2\mbE[|X_{n,g(i)}X_{n,g(i)+1}|]\}\nonumber\\
&\leq&\frac{1}{B_n^2}\sum_{i=1}^{h-1}\{\mbE[X_{n,g(i)}^2]+\mbE[X_{n,g(i)}^2+X_{n,g(i)-1}^2]+\mbE[X_{n,g(i)}^2+X_{n,g(i)+1}^2]\}\nonumber\\
&\leq&\frac{1}{B_n^2}\sum_{i=1}^{h-1}\{3\mbE[X_{n,g(i)}^2]+\mbE[X_{n,g(i)-1}^2]+\mbE[X_{n,g(i)+1}^2]\}\leq3\sum_{i=1}^{h-1}\beta_{n,g(i)}\rightarrow0,\label{eq3.31}\\
	\left|\sum_{i=1}^{h-1}\ud_{n,g(i)}\right|&\leq&\frac{1}{B_n^2}\sum_{i=1}^{h-1}\{\mbE[X_{n,g(i)}^2]+2\mbE[|X_{n,g(i)}X_{n,g(i)-1}|]+2\mbE[|X_{n,g(i)}X_{n,g(i)+1}|]\}\nonumber\\
&\leq&3\sum_{i=1}^{h-1}\beta_{n,g(i)}\rightarrow0.
\ees
Now denote $g(h)=k_n+1$ and
\bess
H_j&=&\{k\in\mathbb{N}:g(j-1)<k<g(j)\},\enspace j=1,\cdots,h,\\
Y_{n,i}&=&\sum_{j\in H_i}X_{n,j},\enspace i=1,\cdots,h,\\
\tilde{B}_{n,K}^2&=&\sum_{i=1}^K \mbE[Y_{n,i}^2],\tilde{b}_{n,K}^2=\sum_{i=1}^K\mbe[Y_{n,i}^2],\enspace K=1,\cdots,h,\enspace \tilde{B}_n^2=\tilde{B}_{n,h},\tilde{b}_n^2=\tilde{b}_{n,h}^2,\\
B_{n,M}^2&=&\mbE[(\sum_{k=1}^M X_{n,k})^2],b_{n,M}^2=\mbe[(\sum_{k=1}^M X_{n,k})^2\enspace M=1,\cdots,k_n,\enspace b_n^2=b_{n,k_n}^2.
\eess

Next, we need to verify that the independent sequence $\{Y_{n,i};i=1,\cdots,h\}$ satisfies the conditions of Lemma \ref{lemma1d}. For condition (\ref{var}), it follows that for any $K=1,\cdots,h$,
\bes
&&\frac{\tilde{B}_{n,K}^2}{B_n^2}=\frac{1}{B_n^2}\sum_{i=1}^K\mbE\left[\left(\sum_{j\in H_i}X_{n,j}\right)^2\right]=\frac{1}{B_n^2}\mbE\left[\sum_{i=1}^K\left(\sum_{j\in H_i}X_{n,j}\right)^2\right]\nonumber\\
&=&\frac{1}{B_n^2}\mbE\left[\left(\sum_{k=1}^{g(K)-1}X_{n,k}\right)^2-\sum_{i=1}^{K-1}\{X_{n,g(i)}^2+2X_{n,g(i)}X_{n,g(i)-1}+2X_{n,g(i)}X_{n,g(i)+1}\}-2\sum_{(i,j)\in\Lambda_K}X_{n,i}X_{n,j}\right]\nonumber\\
&\leq&\frac{B_{n,g(K)-1}^2}{B_n^2}-\sum_{i=1}^{K-1}\ld_{n,g(i)}+\frac{2}{B_n^2}\mbE\left[-\sum_{(i,j)\in\Lambda_K}X_{n,i}X_{n,j}\right],\\
&&\frac{\tilde{B}_{n,K}^2}{B_n^2}\geq\frac{B_{n,g(K)-1}^2}{B_n^2}-\sum_{i=1}^{K-1}\ud_{n,g(i)}-\frac{2}{B_n^2}\mbE\left[\sum_{(i,j)\in\Lambda_K}X_{n,i}X_{n,j}\right],\label{eq3.33}
\ees
where $\Lambda_K=\{(i,j):j-i>1,(i,j)\notin H_l\times H_l,\forall l=1\cdots,K,i,j=1,\cdots,g(K)-1\}$. Denote $\Lambda_{K,j}=\{i:(i,j)\in\Lambda_K,i=1,\cdots,g(K)-1\}$. Note that for $j\in \tilde{H}_l:=H_l\cup\{g(l-1)\},l=2,\cdots,K$, $\Lambda_{K,j}=\{1,\cdots,(j-2)\wedge g(l-1)\}$, for $j\in H_1$, $\Lambda_j=\emptyset$ and for $i\in \Lambda_{K,j}$, $X_{n,j}$ is independent of $X_{n,i}$, thus
\bes
	\mbE\left[\pm\sum_{(i,j)\in\Lambda_K} X_{n,i}X_{n,j}\right]&=&\mbE\left[\sum_{j=1}^{g(K)-1}\sum_{i\in\Lambda_{K,j}} \pm X_{n,i}X_{n,j}\right]\leq \sum_{j=1}^{g(K)-1} \mbE\left[\pm\sum_{i\in\Lambda_{K,j}} X_{n,i}X_{n,j}\right]\nonumber\\
&=&\sum_{j=1}^{g(K)-1}\mbE\left[\mbE[xX_{n,j}]_{x=\pm \sum_{i\in\Lambda_{K,j}}X_{n,i}}\right]\nonumber\\
&\leq&\sum_{j=1}^{g(K)-1}\mbE\left[\left|\sum_{i\in\Lambda_{K,j}}X_{n,i}\right|\right]\left(|\mbE[X_{n,j}]|+|\mbe[X_{n,j}]|\right)\nonumber\\
&\leq&\mbE\left[\max_{k\leq k_n}\left|\sum_{i=1}^kX_{n,i}\right|\right]\sum_{j=1}^{k_n}\left(|\mbE[X_{n,j}]|+|\mbe[X_{n,j}]|\right),\label{eq3.35}
\ees
where we define that the summation over empty set is zero and by Lemma \ref{lemma4.1},
\bes
	\mbE\left[\max_{k\leq k_n}\left|\sum_{i=1}^kX_{n,i}\right|\right]&\leq&\mbE\left[\max_{k\leq k_n}\left|\sum_{i=1}^kX_{n,i}\right|^2\right]^{1/2}\nonumber\\
&\leq& C\left\{\sum_{k=1}^{k_n}\mbE[|X_{n,k}|^2]+\left(\sum_{k=1}^{k_n}\left[|\mbE[X_{n,k}]|+|\mbe[X_{n,k}]|\right]\right)^2\right\}^{1/2}\nonumber\\
&=&O(B_n).\label{eq3.36}
\ees
Combining (\ref{eq3.31})-(\ref{eq3.36}) we have
\bes
	\frac{\tilde{B}_{n,K}^2}{B_n^2}-\frac{B_{n,g(K)-1}^2}{B_n^2}\rightarrow0.\label{BB}
\ees
Let $K=h$, we have $\frac{\tilde{B}_n^2}{B_n^2}\rightarrow1$,
and similarly
\bes
	\frac{\tilde{b}_{n,K}^2}{B_n^2}-\frac{b_{n,g(K)-1}^2}{B_n^2}\rightarrow 0.\label{bb}
\ees
When $\liminf\frac{\tilde{B}_{n,K}^2}{\tilde{B}_n^2}>0$, it follows that 
\bes
	\frac{\tilde{B}_{n,K}^2}{\tilde{B}_n^2}-\frac{B_{n,g(K)-1}^2}{B_n^2}=\frac{\tilde{B}_{n,K}^2}{B_n^2}\cdot\frac{B_n^2}{\tilde{B}_n^2}-\frac{B_{n,g(K)-1}^2}{B_n^2}\rightarrow0,\label{diff}
\ees
thus $\liminf\frac{B_{n,g(K)-1}^2}{B_n^2}>0$ as $h>K\rightarrow\infty$, and by condition (\ref{eq3.29}),
\bes
\frac{b_{n,g(K)-1}^2}{B_{n,g(K)-1}^2}\rightarrow r.	\label{cc}
\ees
Note that
\bess
	\left(\frac{\tilde{b}_{n,K}^2}{\tilde{B}_{n,K}^2}-r\right)\cdot\frac{\tilde{B}_{n,K}^2}{\tilde{B}_n^2}+\left(r\cdot\frac{\tilde{B}_{n,K}^2}{\tilde{B}_n^2}-\frac{b_{n,g(K)-1}^2}{B_{n,g(K)-1}^2}\cdot\frac{B_{n,g(K)-1}^2}{B_n^2}\right)=\frac{\tilde{b}_{n,K}^2}{B_n^2}-\frac{b_{n,g(K)-1}^2}{B_n^2}\rightarrow 0,
\eess
it follows from (\ref{diff}) and (\ref{cc}) that
\bess
\frac{\tilde{b}_{n,K}^2}{\tilde{B}_{n,K}^2}\rightarrow r
\eess
with $h>K\rightarrow\infty$ as long as $\liminf\frac{\tilde{B}_{n,K}^2}{\tilde{B}_n^2}>0$. Hence condition (\ref{var}) is satisfied.
For the Lindeberg condition (\ref{Lin}), we have
\bess
	&&\frac{1}{B_n^2}\sum_{i=1}^{h}\mbE[(Y_{n,i}^2-\epsilon^2B_n^2)^+]\\
&=&\frac{1}{B_n^2}\sum_{i=1}^{h}\mbE\left[\left(\left(\sum_{j\in H_i}X_{n,j}\right)^2-\epsilon^2B_n^2\right)^+\right]\leq\frac{1}{B_n^2}\sum_{i=1}^{h}\mbE\left[\left(p_n\sum_{j\in H_i}X_{n,j}^2-\epsilon^2B_n^2\right)^+\right]\\
&\leq& \frac{p_n}{B_n^2}\sum_{i=1}^{h}\mbE\left[\left(\sum_{j\in H_i}(X_{n,j}^2-\frac{\epsilon^2B_n^2}{p_n^2})\right)^+\right]\leq\frac{p_n}{B_n^2}\sum_{i=1}^{h}\sum_{j\in H_i}\mbE[(X_{n,j}^2-\frac{\epsilon^2B_n^2}{p_n^2})^+]\rightarrow0,
\eess
where the last inequality is due to $(a+b)^+\leq a^++b^+$ and the sub-linearity of $\mbE$. Since $\frac{\tilde{B}_n^2}{B_n^2}\rightarrow 1$, it is obvious that
\bess
	\frac{1}{\tilde{B}_n^2}\sum_{i=1}^{h}\mbE[(Y_{n,i}^2-\epsilon^2\tilde{B}_n^2)^+]&=&\frac{\tilde{B}_n^2}{B_n^2}\cdot\frac{1}{\tilde{B}_n^2}\sum_{i=1}^{h}\mbE\left[\left(Y_{n,i}^2-\epsilon^2B_n^2\cdot\frac{\tilde{B}_n^2}{B_n^2}\right)^+\right]\\
&\leq&\frac{2}{B_n^2}\sum_{i=1}^{h}\mbE[(Y_{n,i}^2-\epsilon^2B_n^2/2)^+]\rightarrow0.
\eess
 For (\ref{mean}), note that 
\bess
\sum_{j\in H_i}\mbe[X_{n,j}]\leq\mbe\left[\sum_{j\in H_i}X_{n,j}\right]\leq\mbE\left[\sum_{j\in H_i}X_{n,j}\right]\leq\sum_{j\in H_i}\mbE[X_{n,j}],
\eess
and denote
\bess
I_1=\left\{i:\mbE\left[\sum_{j\in H_i}X_{n,j}\right]\geq0\right\},\enspace I_2=\left\{i:\mbE\left[\sum_{j\in H_i}X_{n,j}\right]<0\right\},
\eess
it follows that
\bess
\sum_{i=1}^h\left|\mbE\left[\sum_{j\in H_i}X_{n,j}\right]\right|&=&\sum_{i\in I_1}\mbE\left[\sum_{j\in H_i}X_{n,j}\right]-\sum_{i\in I_2}\mbE\left[\sum_{j\in H_i}X_{n,j}\right]\\
&\leq& \sum_{i\in I_1}\sum_{j\in H_i}\mbE[X_{n,j}]-\sum_{i\in I_2}\sum_{j\in H_i}\mbe[X_{n,j}]\\
&\leq& \sum_{i\in I_1}\sum_{j\in H_i}|\mbE[X_{n,j}]|+\sum_{i\in I_2}\sum_{j\in H_i}|\mbe[X_{n,j}]|\\
&\leq&\sum_{k=1}^{k_n}\{|\mbE[X_{n,j}]|+|\mbe[X_{n,j}]|\}=o(B_n)=o(\tilde{B}_n).
\eess
The lower expectation part can be inferred by analogy, thus the condition (\ref{mean}) is satisfied for $\{Y_{n,i};i=1,\cdots,h\}$. 
Hence it follows from Lemma \ref{lemma1d} that for any continuous function $\varphi$ with $|\varphi(x)|\leq Cx^2$,
\bes
\lim_{n\rightarrow\infty}\mbE\left[\varphi\left(\frac{\sum_{i=1}^h Y_{n,i}}{\tilde{B}_n}\right)\right]=\tE[\varphi(\xi)],
\ees
where $\xi\sim N(0,[r,1])$. On the other hand, $\{X_{n,g(i)};i=1,\cdots,h-1\}$ is an independent sequence and by the proof of Theorem 2.1 of Zhang\cite{zhang23},
\bess
\lim_{n\rightarrow\infty}\mbE\left[\left(\sum_{i=1}^{h-1}X_{n,g(i)}\right)^2\right]=\lim_{n\rightarrow\infty}\frac{1}{B_n^2}\sum_{i=1}^{h-1}\mbE[X_{n,g(i)}^2]\leq\lim_{n\rightarrow\infty}\sum_{i=1}^{h-1}\beta_{n,g(i)}=0.
\eess
Thus $\sum_{i=1}^{h-1}X_{n,g(i)}/B_n\overset{\V}{\rightarrow}0$ and by Slutsky's Theorem (c.f. Lemma4.2 of Zhang\cite{zhang15}), the theorem for $m=1$ is proved for bounded continuous function $\varphi$.
\par Now condider the general case. Let $k_n'=[\frac{k_n}{m}]+1$ and
\bes
	Z_{n,k}=\sum_{i=1}^mX_{n,m(k-1)+i},\enspace k=1\cdots,k_n'-1,\quad Z_{n,k_n'}=\sum_{k=m(k_n'-1)+1}^{k_n}X_{n,k}.\label{eq2.43}
\ees
It is obvious that $\{Z_{n,k};k=1\cdots,k_n'\}$ is 1-dependent and it is easy to verify that conditions (\ref{eq3.26})-(\ref{eq3.29}) hold for $\{Z_{n,k};k=1,\cdots,k_n'\}$.
\par Suppose that (\ref{eqpp}) holds, with the same argument of Zhang\cite{zhang23}, it is sufficient to show that
\bes
	\lim_{N\rightarrow\infty}\limsup_{n\rightarrow\infty}\mbE\left[\left(\left|\frac{1}{B_n}\sum_{k=1}^{k_n}X_{n,k}\right|^p-N\right)^+\right]=0.
\ees
Without loss of generality, we can assume that $B_n=1$ by setting $X'_{n,k}=X_{n,k}/B_n$. Let $\hat{X}_{n,k}=(-1)\vee\Xnk\wedge1$ and $\tilde{X}_{n,k}=\Xnk-\hat{X}_{n,k}$. It is obvious that
\bes
	&&\sum_{k=1}^{k_n}\mbE[|\tilde{X}_{n,k}|]=\sum_{k=1}^{k_n}\mbE\left[(|X_{n,k}|-1)^+\right]\leq2\sum_{k=1}^{k_n}\mbE\left[(|X_{n,k}^2-1/2|)^+\right]\rightarrow 0,\\
	&&\sum_{k=1}^{k_n}\left\{|\mbE[\hat{X}_{n,k}]|+|\mbe[\hat{X}_{n,k}]|\right\}\rightarrow 0,\label{eq3.46}\\
	&&\sum_{k=1}^{k_n}\mbE[|\hat{X}_{n,k}|^q]\leq\sum_{k=1}^{k_n}\mbE[X_{n,k}^2]=O(1),\forall q\geq2.\label{eq3.47}
\ees
For $m$-dependent random variables, it follows from  Lemma \ref{lemma4.1} for $q\geq2$ that,
\bess
\mbE\left[\left|\sum_{k=1}^{k_n}\hat{X}_{n,k}\right|^q\right]
\leq C_{m,q}\left\{\sum_{k =1}^{k_n}\mbE[|\hat{X}_{n,k}|^q]+\left(\sum_{k=1}^{k_n}\mbE[\hat{X}_{n,k}^2]\right)^{q/2}+\left(\sum_{k=1}^{k_n} \left(|\mbE[\hat{X}_{n,k}]|+|\mbe[\hat{X}_{n,k}]|\right)\right)^q\right\}.
\eess
Thus we have
\bess
\mbE\left[\left|\sum_{k=1}^{k_n}\hat{X}_{n,k}\right|^q\right]=O(1),
\eess
by (\ref{eq3.46}) and (\ref{eq3.47}). Then we have
\bess
	\lim_{N\rightarrow\infty}\limsup_{n\rightarrow\infty}\mbE\left[\left(\left|\sum_{k=1}^{k_n}\hat{X}_{n,k}\right|^p-N\right)^+\right]\\
\leq\lim_{N\rightarrow\infty}\limsup_{n\rightarrow\infty}N^{-1}\mbE\left[\left|\sum_{k=1}^{k_n}\hat{X}_{n,k}\right|^{2p}\right]=0.
\eess
And similarly for $\tilde{X}_{n,k}$, by the Rosenthal-type inequality again, we have
\bess
	\mbE\left[\left|\sum_{k=1}^{k_n}\tilde{X}_{n,k}\right|^p\right]\leq C_{m,p}\left\{\sum_{k=1}^{k_n}\mbE[(|X_{n,k}|^p-1)^+]+\left(\sum_{k=1}^{k_n}\mbE[(X_{n,k}^2-1)^+]\right)^{p/2}+\left(\sum_{k=1}^{k_n}\mbE[(|\Xnk|-1)^+]\right)^p\right\}\rightarrow0.
\eess
The proof is completed.\qedsymbol
\par The following is a direct corollary of Theorem \ref{th3.3}, where we consider a sequence of random variables instead of arrays of random variables.
\begin{coro}
	Let $\{X_k;k=1,\cdots,n\}$ be a sequence of $m$-dependent random variables in the sub-linear expectation space $\sles$ with $\mbE[X_k^2]<\infty$ for $k=1,\cdots,n$. Denote $S_n=\sum_{k=1}^n X_k$, $B_n^2=\mbE[S_n^2] $ and $b_n^2=\mbe[S_n^2]$. Assume that
\bes
	&&\frac{1}{B_n^2}\sum_{k=1}^n \mbE[(X_k^2-\epsilon B_n^2)^+]\rightarrow 0,\enspace \forall \epsilon>0,\label{eq3.48}\\
 	&&\frac{1}{B_n}\sum_{k=1}^n \{|\mbE[X_k]|+|\mbe[X_k]|\}\rightarrow0,\label{eq3.49}\\
	&&\frac{1}{B_n^2}\sum_{k=1}^n \mbE[X_k^2]=O(1),\label{eq3.50}
\ees 
and there exists $r\in[0,1]$ such that
\bes
	b_n^2/B_n^2\rightarrow r.\label{eq3.51}
\ees
Then for any continuous function $\varphi\in C(\mathbb{R})$ with $|\varphi(x)|\leq Cx^2$, we have
\bes
	\lim_{n\rightarrow\infty}\mbE\left[\varphi\left(\frac{S_n}{B_n}\right)\right]=\tE[\varphi(\xi)],\label{eq3.53}
\ees
where $\xi\sim N(0,[r,1])$. Further, when $p>2$, (\ref{eq3.53}) holds for any continuous function $\varphi\in C(\mathbb{R})$ with $|\varphi(x)|\leq C|x|^p$ if (\ref{eq3.48}) is replaced by the condition that 
\bes
	\frac{1}{B_n^p}\sum_{k=1}^{k_n}\mbE[|X_{k}|^p]\rightarrow0.
\ees
\label{coro2.1}
\end{coro}
\begin{remark}
	In contrast to Theorem 3.4 in Zhang\cite{zhang23}, if $\{X_k;k=1\cdots,n\}$ is a sequence of independent random variables, under condition (\ref{eq3.49}), we have
\bess
	B_n^2=\mbE\left[\left(\sum_{k=1}^nX_k\right)^2\right]=\sum_{k=1}^n\mbE[X_k^2]+o(B_n^2),
\eess
thus (\ref{eq3.50}) is implied by (\ref{eq3.49}) and can be removed under independent situation and Corollary \ref{coro2.1} coincides with Theorem 3.4 of Zhang\cite{zhang23}. 
\end{remark}

If we consider a special case without mean uncertainty, (\ref{eq3.49}) is satisfied obviously. We can rewrite Corollary \ref{coro2.1} into the following form, which improves Theorem 3.1 of Li\cite{li15}.
\begin{coro}
Let $\{X_k;k=1,\cdots,n\}$ be a sequence of $m$-dependent random variables in the sub-linear expectation space $\sles$ with $\mbE[X_k^2]<\infty$ for $k=1,\cdots,n$. Denote $S_n=\sum_{k=1}^n X_k$, $B_n^2=\mbE[S_n^2] $ and $b_n^2=\mbe[S_n^2]$. Assume that
\bes
	&& \mbE[(X_n^2-\epsilon n)^+]\rightarrow 0,\enspace \forall \epsilon>0,\label{eq2.35}\\
 	&&\mbE[X_k]=\mbe[X_k]=0,
\ees
 and there exist constants $0\leq\lsigma^2\leq\usigma^2<\infty$ such that
\bes
	\frac{B_n^2}{n}\rightarrow\usigma^2,\enspace \frac{b_n^2}{n}\rightarrow\lsigma^2.\label{eq2.37}
\ees
Then for any continuous function $\varphi\in C(\mathbb{R})$ with $|\varphi(x)|\leq Cx^2$, we have
\bess
	\lim_{n\rightarrow\infty}\mbE\left[\varphi\left(\frac{S_n}{\sqrt{n}}\right)\right]=\tE[\varphi(\xi)],
\eess
where $\xi\sim N(0,[\lsigma^2,\usigma^2])$.
\label{coro2.2}
\end{coro}

Now consider that $\{X_k;k=1,\cdots,n\}$ is a stationary sequence, we can further weaken the conditions (\ref{eq2.35})-(\ref{eq2.37}).
\begin{coro}
Let $\{X_k;k=1,\cdots,n\}$ be a stationary sequence of $m$-dependent random variables in the sub-linear expectation space $\sles$ with $\mbE[X_k^2]<\infty$ for $k=1,\cdots,n$. Denote $S_n=\sum_{k=1}^n X_k$, $B_n^2=\mbE[S_n^2] $ and $b_n^2=\mbe[S_n^2]$. Assume that $\mbE[X_1]=\mbe[X_1]=0$, and
\bes
\mbE[(X_1^2-c)^+]\rightarrow 0,\enspace as\enspace c\rightarrow\infty.\label{eq2.38}
\ees
Then (\ref{eq2.35}) and (\ref{eq2.37}) hold and thus for any continuous function $\varphi\in C(\mathbb{R})$ with $|\varphi(x)|\leq Cx^2$, we have
\bess
	\lim_{n\rightarrow\infty}\mbE\left[\varphi\left(\frac{S_n}{\sqrt{n}}\right)\right]=\tE[\varphi(\xi)],
\eess
where $\xi\sim N(0,[\lsigma^2,\usigma^2])$ for some $0\leq\lsigma^2\leq\usigma^2<\infty$.
\label{coro2.3}
\end{coro}
\proof  It is obvious that
\bes
	 \mbE[S_n]=\mbe[S_n]=0
\ees
by Proposition 1.3.7 of Peng\cite{Peng19}, and
\bess
	B_n^2=\mbE[S_n^2]&=&\mbE\left[\sum_{k=1}^nX_k^2+2\sum_{1\leq i<j\leq n}X_iX_j\right]\\
&\leq&\sum_{k=1}^n\mbE[X_k^2]+2\sum_{1\leq i<j\leq n}\mbE[X_iX_j]\\
&\leq&\sum_{k=1}^n\mbE[X_k^2]+2\sum_{i=1}^n\sum_{j=i+1}^{n\wedge (i+m)}\mbE[X_iX_j]\leq nC.
\eess
So there is a constant $0\leq \usigma^2<\infty$ and a sequence $0<p_n\uparrow\infty$ such that
\bes
	\frac{\mbE[S_{p_n}^2]}{p_n}\rightarrow \usigma^2.
\ees
Without loss of generality, we can assume that $p_n=o(n)$. Otherwise, we can set $k_n=\sup\{k:p_k\leq\sqrt{n}\}$ and redefine $p_k$ by $p_{k_n}$.We devide $S_n$ into three parts and denote them by $A_{n,j},j=1,2,3$, respectively, then we have
\bes
	&&\lim_{n\rightarrow\infty}\frac{1}{n}\mbE[A_{n,1}^2]=\lim_{n\rightarrow\infty}\frac{1}{n}\mbE\left[\left(\sum_{i=0}^{[\frac{n}{p_n+m}]-1}\sum_{j=1}^{p_n}X_{i(p_n+m)+j}\right)^2\right]\nonumber\\
&&=\lim_{n\rightarrow\infty}\frac{1}{n}\left[\frac{n}{p_n+m}\right]\mbE[S_{p_n}^2]=\lim_{n\rightarrow\infty}\frac{p_n}{n}\left[\frac{n}{p_n+m}\right]\frac{\mbE[S_{p_n}^2]}{p_n}=\usigma^2,\label{eq2.41}\\
&&\lim_{n\rightarrow\infty}\frac{1}{n}\mbE[A_{n,2}^2]=\lim_{n\rightarrow\infty}\frac{1}{n}\mbE\left[\left(\sum_{i=0}^{[\frac{n}{p_n+m}]-1}\sum_{j=p_n+1}^{p_n+m}X_{i(p_n+m)+j}\right)^2\right]\nonumber\\
&&=\lim_{n\rightarrow\infty}\frac{m}{n}\left[\frac{n}{p_n+m}\right]\frac{\mbE[S_{m}^2]}{m}=0,
\ees
and 
\bes
\limsup_{n\rightarrow\infty}\frac{1}{n}\mbE[A_{n,3}^2]=\limsup_{n\rightarrow\infty}\frac{1}{n}\mbE\left[\left(\sum_{j=[\frac{n}{p_n+m}](p_n+m)}^{n}X_{j}\right)^2\right]\leq \limsup_{n\rightarrow\infty}\frac{p_n+m}{n}C=0.
\ees
Note that 
\bes
	\left|\Vert A_{n,1}+A_{n,2}+A_{n,3}\Vert_2-\Vert A_{n,1}\Vert_2\right|\leq\Vert A_{n,2}+A_{n,3}\Vert_2\leq\Vert A_{n,2}\Vert_2+\Vert A_{n,3}\Vert_2,\label{eq2.44}
\ees
where $\Vert X\Vert_2=\mbE[X^2]^{1/2}$,it follows from (\ref{eq2.41})-(\ref{eq2.44}) that
\bess
	\frac{\mbE[S_n^2]}{n}\rightarrow\usigma^2,
\eess
and similarly we can obtain that
\bess
	\frac{\mbe[S_n^2]}{n}\rightarrow\lsigma^2.
\eess
Hence (\ref{eq2.37}) holds and (\ref{eq2.35}) is trivial from (\ref{eq2.37}) and (\ref{eq2.38}). \qedsymbol
\begin{remark}
	In particular, if we consider the i.i.d. case, Corollary \ref{coro2.1}-\ref{coro2.3} coincide with Peng's central limit theorem(c.f. Theorem 2.4.4 of Peng\cite{Peng19}).
\end{remark}

\section{Central Limit Theorems for truncated conditions}
From Lemma \ref{lemma1d}, we have the following central limit theorem for one-dimensional random variables with truncated conditions. 
\begin{theorem}
	Let $\{X_{n,k};k=1,\cdots,k_n\}$ be an array of  independent random variables in the sub-linear expectation space $\sles$. Assume that there exists a constant $\tau>0$ such that 
\bes
	\sum_{k=1}^{k_n}\V(|X_{n,k}|>\epsilon )\rightarrow0\enspace\forall\epsilon>0,\label{eq2.6}\\
	\frac{1}{B_n}\sum_{k=1}^{k_n}\left\{|\mbE[X_{n,k}^{(\tau)}]|+|\mbe[X_{n,k}^{(\tau)}]|\right\}\rightarrow0,\label{eq2.7}
\ees
where $B_n^2=\sum_{k=1}^{k_n}\mbE[(X_{n,k}^{(\tau)})^2]$ and assume that there is a constant $r\in[0,1]$ such that
\bes
	\sum_{k=1}^{m}\mbe[(X_{n,k}^{(\tau)})^2]/\sum_{k=1}^{m}\mbE[(X_{n,k}^{(\tau)})^2]\rightarrow r,\label{eq2.8}
\ees
as long as $k_n>m\rightarrow \infty$ with $\liminf\sum_{k=1}^{m}\mbE[(X_{n,k}^{(\tau)})^2]/B_n^2>0.$
Then for any bounded continuous function $\varphi\in C_b(\mathbb{R})$, we have
\bes
	\lim_{n\rightarrow\infty}\mbE\left[\varphi\left(\frac{\sum_{k=1}^{k_n}X_{n,k}}{B_n}\right)\right]=\tE[\varphi(\xi)],\label{eq2.9}
\ees
where $\xi\sim N(0,[r,1])$. 
\label{th2.2}
\end{theorem}
\proof From (\ref{eq2.6}), it follows that 
\bess
	\V(\Xnk\neq \Xnk^{(\tau)}\enspace for \enspace some \enspace k)=\V(\max_{1\leq k\leq k_n}|\Xnk|>\tau)\leq\sum_{k=1}^{k_n}\V(|\Xnk|>\tau)\rightarrow 0.
\eess
For any bounded function $\varphi$, we have
\bess
	\mbE\left[\left|\varphi\left(\frac{\sum_{k=1}^{k_n}\Xnk}{B_n}\right)-\varphi\left(\frac{\sum_{k=1}^{k_n}\Xnk^{(\tau)}}{B_n}\right)\right|\right]\leq 2\sup_x|\varphi(x)|\V(\Xnk\neq \Xnk^{(\tau)}\enspace for \enspace some \enspace k)\rightarrow 0.
\eess
Hence it is sufficient to verify the Linderberg condition (\ref{Lin}) for $\{\Xnk^{(\tau)}\}$. For any $\epsilon>0$,
\bess
\sum_{k=1}^{k_n}\mbE[((\frac{\Xnk^{(\tau)}}{B_n})^2-\epsilon)^+]\leq (\frac{\tau}{B_n})^2\sum_{k=1}^{k_n}\V(|X_{n,k}|>\epsilon)\rightarrow0,
\eess
thus (\ref{eq2.9}) follows from Lemma \ref{lemma1d}.\qedsymbol
\par Taking advantage of Theorem \ref{th2.2}, we give the central limit theorem for $m$-dependent random variables without Lindeberg condition in sub-linear expectation space.  
\begin{theorem}
Let $\{X_{n,k};k=1,\cdots,k_n\}$ be an array of  $m$-dependent random variables in the sub-linear expectation space $\sles$. Assume that there exist constants $r\in[0,1]$ and $\tau>0$ such that 
\bes
	\sum_{k=1}^{k_n}\V(|X_{n,k}|>\epsilon)\rightarrow0\enspace\forall\epsilon>0,\label{eq3.1}\\
	\frac{1}{B_n}\sum_{k=1}^{k_n}\left\{|\mbE[X_{n,k}^{(\tau)}]|+|\mbe[X_{n,k}^{(\tau)}]|\right\}\rightarrow0,\label{eq3.2}\\
	\frac{1}{B_n^2}\sum_{k=1}^{k_n}\mbE[(X_{n,k}^{(\tau)})^2]=O(1),\label{eq3.3}
\ees
and
\bes
	\sum_{k=1}^{M}\mbe[(X_{n,k}^{(\tau)})^2]/\sum_{k=1}^{M}\mbE[(X_{n,k}^{(\tau)})^2]\rightarrow r,\label{eq3.4}
\ees
as long as $k_n>M\rightarrow \infty$ with $\liminf\sum_{k=1}^{M}\mbE[(X_{n,k}^{(\tau)})^2]/B_n^2>0,$ where $B_n^2=\sum_{k=1}^{k_n}\mbE[(X_{n,k}^{(\tau)})^2]$.
Then for any bounded continuous function $\varphi\in C_b(\mathbb{R})$, we have
\bes
	\lim_{n\rightarrow\infty}\mbE\left[\varphi\left(\frac{\sum_{k=1}^{k_n}X_{n,k}}{B_n}\right)\right]=\tE[\varphi(\xi)],\label{eq3.5}
\ees
where $\xi\sim N(0,[r,1])$.
\label{th3.1}
\end{theorem}
\proof The proof is similar to Theorem \ref{th3.3}, so we only give a sketch for $m=1$. From (\ref{eq3.1}), there exists a non-decreasing sequence $p_n\uparrow\infty,p_n=o(n)$ satisfying
\bes
	\sum_{k=1}^{k_n}\V(|X_{n,k}|>\frac{\epsilon}{p_n})\rightarrow0,\enspace\forall\epsilon>0.\label{eq3.5}
\ees
Let $\tau'=\frac{\tau}{p_n}$, it is easily verified that all the conditions (\ref{eq3.2})-(\ref{eq3.4}) hold for $\{X_{n,k}^{(\tau')}\}$. Denote 
\bess
\beta_{n,k}^{\tau}&=&\{\mbE[(X_{n,k-1}^{(\tau)})^2]+\mbE[(\Xnk^{(\tau)})^2]+\mbE[(X_{n,k+1}^{(\tau)})^2]\}/B_n^2,\\
\ld_{n,k}^\tau&=&\{\mbe[(\Xnk^{(\tau)})^2]+2\mbe[X_{n,k}^{(\tau)}X_{n,k-1]}^{(\tau)}]+2\mbe[X_{n,k}^{(\tau)}X_{n,k+1]}^{(\tau)}]\}/B_n^2,\\
\ud_{n,k}^\tau&=&\{\mbE[(\Xnk^{(\tau)})^2]+2\mbE[X_{n,k}^{(\tau)}X_{n,k-1]}^{(\tau)}]+2\mbE[X_{n,k}^{(\tau)}X_{n,k+1]}^{(\tau)}]\}/B_n^2,
\eess
and define $P_i,g(i)$ and $H_i$ in the same way. Let $Y_{n,i}=\sum_{j\in H_i}X_{n,j}$, we can compare $Y_{n,i}^{(\tau)}$ with $\sum_{j\in H_i}X_{n,j}^{(\tau')}$ and show that all the conditions in Theorem \ref{th2.2} are satisfied for $Y_{n,i}^{(\tau)}$. Hence it follows from Theorem \ref{th2.2} that for any bounded continuous function $\varphi$,
\bes
\lim_{n\rightarrow\infty}\mbE\left[\varphi\left(\frac{\sum_{i=1}^h Y_{n,i}}{B_n}\right)\right]=\tE[\varphi(\xi)],
\ees
where $\xi\sim N(0,[r,1])$. After the same procedure, we can show that 
\bess
	\sum_{i=1}^{h-1}X_{n,g(i)}/B_n\overset{\V}{\rightarrow}0.
\eess
Hence we obtain the conclusion for $m=1$. Take $Z_{n,k}$ as (\ref{eq2.43}), the $m$-dependent case can be inferred immediately .\qedsymbol

\bibliographystyle{plain}
\bibliography{m-CLT}

\end{document}